\documentclass[12pt,a4paper]{article}
\usepackage[top=3cm, bottom=3cm, left=3cm, right=3cm]{geometry}
\usepackage[latin1]{inputenc}
\usepackage{amsmath}
\usepackage{amsthm}
\usepackage{amsfonts}
\usepackage{amssymb}
\usepackage{graphicx}

\usepackage{amsmath,amsthm,amssymb,graphicx, multicol, array}
\usepackage{enumerate}
\usepackage{hyperref}

\newtheorem{thm}{Theorem}[section]
\newtheorem{lem}{Lemma}[section]

\newcommand{\Z}{\mathbb{Z}}

\newcommand{\R}{\mathbb{R}}

\title{The Error Term of the Summatory Euler Phi Function}
\date{}
\author{N. A. Carella}
\begin{document}
\maketitle

\tableofcontents

\section{Introduction} \label{sec1}
The error term is defined by $E(x)=\sum _{n \leq  x} \varphi (n)-3\pi ^{-2}x$. The earliest estimate for the error term $E(x)=O\left(x^{1+\epsilon}\right)$, with $\epsilon >0$, was computed by Mertens, followed by $E(x)=O(x \log  x)$ computed by Dirichlet. The current error term $E(x)=O\left(x(\log
x)^{2/3}(\log \log x)^{4/3}\right)$ in the mathematical literature is attributed to Walfisz, confer \cite[p.\ 68]{MV07}, \cite[p.\ 102]{SH83}, \cite[p.\ 47]{TG15},  et alii. For a large real number $x\in \mathbb{R}$, the explicit formula
\begin{equation} \label{1000}
\sum_{n\leq x} \varphi(n)=\frac{1}{6}+\frac{\delta(x)}{2}\varphi(x)+\frac{3}{\pi^2}x^{2}+\sum_{\rho}\frac{\zeta (\rho -1)}{\rho \zeta '(\rho )}x^{\rho }+
\sum_{n\geq1}\frac{\zeta (-2n-1)}{-2n \zeta '(-2n)}x^{-2n},
\end{equation}
where
\begin{equation}
\delta(n) =
\left\{
\begin{array}{lr}
1     &x \in \mathbb{N},\\
0      &x \notin \mathbb{N},
\end{array}
\right .
\end{equation}

shows that there is a sharper unconditional error term 
\begin{equation}
E(x)=\delta (x) \varphi(x)/2+O\left(x e^{-c\sqrt{\log x}}\right). 
\end{equation}
The proof of (1) is based on standard analytical techniques and the Perron summation formula, see \cite[p.\ 138]{MV07}, \cite[p.\ 217]{TG15}, and similar references. The result in Theorem \ref{thm1.1} provides a different and independent proof of this sharper unconditional error term.
\\

\begin{thm} \label{thm1.1}  For a large number $x\geq 1$, the average order for the Euler totient function $\varphi (n)$ has the asymptotic formula
\begin{equation}
\sum_{n\leq x} \varphi(n)=\frac{3}{\pi^2}x^{2}+O(x),
\end{equation}
unconditionally.
\end{thm}

\begin{thm} \label{thm1.2}  Assume the RH. For a large number $x\geq 1$, the average order for the Euler totient function $\varphi (n)$ has the asymptotic formula
	\begin{equation}
	\sum_{n\leq x} \varphi(n)=\frac{3}{\pi^2}x^{2}+\frac{\delta(x)}{2}\varphi(x)+O(x^{1/2}\log^2 x).
	\end{equation}
\end{thm}

The next few sections provide the elementary background materials, and the last section has the proof of Theorem \ref{thm1.1}. The proof of Theorem \ref{thm1.2} follows from the explicit formula (\ref{1000}).
\\

\begin{thm} \label{thm1.3}   Let $x \geq 1$ be a large number. Then, the average order of the ratio $1/\varphi(n)$ is as follows.
	
	\begin{equation} \label{90013}	
	\sum_{n \leq x}  \frac{1}{\varphi(n)} =c_0+c_1\log x +O\left (\frac{\log x}{x} \right ),
	\end{equation}	
	where
	$c_0$ and $c_1$ are constants.	
\end{thm}

It is shown that this is the best possible, as was determined by Landau over a century ago, in \cite[p.\ 184]{LE00}. Several different and independent proofs of Theorem \ref{thm1.3} are possible. A proof based on the convolution method is provided in Sections \ref{sec6}.
\section{Results For The Mobius Function}\label{sec2}
The zeta function $\zeta (s)=\sum _{n\geq 1} n^{-s}$ has a pole at $s=1$, so the inverse series $1/\zeta (s)=\sum _{n\geq 1} \mu (n)n^{-s}$ vanishes at $s=1$. For large $x\geq 1$, the associated summatory functions are $\sum _{n\leq x} \mu (n)=o(x)$, and $\sum _{n\leq x} \mu (n)n^{-1}<1$. The true rates of growth and decay respectively of these summatory functions are of considerable interest in number theory.
\\

\begin{lem} \label{lem21001} Let $x\geq 1$ be a large number, and let $\mu(n)$ be the Mobius function. Then
\begin{equation}
\sum_{n\leq x}\mu(n)=O\left (\dfrac{x}{\log^{2} x} \right ).
\end{equation}
\end{lem}
\begin{proof} Refer to \cite[p.\ 182]{MV07}, \cite{IK04} and the literature. 
\end{proof}

A better unconditional result $\sum_{n\leq x} \mu (n)=O\left(xe^{-c(\log x)^{3/5}(\log\log x )^{1/5}}\right)$, with $c>0$ an absolute constant,  is available in the literature, see \cite{FK02}. However, the weaker but sufficient result in Lemma \ref{lem21001}, which has simpler notation, will be used here. 
\\

\begin{lem} \label{lem21040} Let $x\geq 1$ be a large number, let $\mu (n)$ be the Mobius function, and let $
s\in \mathbb{C}$, with $\Re e(s)\geq 1$. Then
\begin{equation}
\sum_{n\leq x} \dfrac{\mu(n)}{n^{s}}=\dfrac{1}{\zeta(s)}+O\left ( \dfrac{1}{x^{s-1} \log^{2}x}  \right ) .
\end{equation}
\end{lem}

\begin{proof} Use Lemma \ref{lem21001} to evaluate the Stieltjes integral

\begin{eqnarray}
\sum_{n\leq x}\frac{\mu(n)}{n^{s}}&=&\sum_{n\geq 1}\frac{\mu(n)}{n^{s}}-\sum_{n\geq x}\frac{\mu(n)}{n^{s}} \nonumber\\
&=& \frac{1}{\zeta(s)}-\int_{x}^{\infty}\frac{1}{t^s} d M(t),
\end{eqnarray}
where \(M(x)=\sum _{n\leq x} \mu (n)\).  
\end{proof}

\section{Finite Fractional Sums} \label{sec3}
Let the symbol \(\{ z \} = z -[ z ]\) denotes the fractional part function. The earliest results for the summatory fractional function are due to Dirichlet and delaVallee Poussin, see \cite{PF10} and \cite[p.\ 259]{MR08} for different proofs. \\

\begin{thm} \label{thm15000} {\normalfont (Dirichlet)}  Let $x \geq 1$ be a number, and let $\{x\}=x-[x]$ be the fractional part function. Then
	\begin{equation}
	\sum_{n \leq x}  \left \{\frac{x}{n} \right \}=(1-\gamma)x + O(x^{1/2}),
	\end{equation}
where \(\gamma =0.5772156649\ldots\) is Euler constant.
\end{thm}

A sharper error term \(O\left(x^{1/3}\right)\) can be deduced from the Voronoi estimate for the divisor problem, the proofs appear in \cite{IM88}, \cite{HM93}, \cite[Section 13.8]{IA03}, \cite[p.\ 201]{CO07}, \cite{CS07}, and other references. 

\begin{lem} \label{lem15003}
	Let $x  \geq 1$ be a large number. Then
	\begin{equation}
	\sum _{n\leq x}\left \{\frac{x}{n} \right \}=(1-\gamma)x+O\left(x^{1/3}\right),
	\end{equation}
	where \(\gamma =0.5772156649\ldots\) is Euler constant.	
\end{lem}

\begin{proof} Express the sum of fractional parts as difference of the harmonic finite sum, and the average number of divisors:
\begin{eqnarray}
\sum _{n\leq x} \left \{\frac{x}{n}\right \}&=&\sum_{n \leq x} \left(\frac{x}{n}-\left [\frac{x}{n}\right ]\right) \nonumber\\
&=&x \sum_{n \leq x} \frac{1}{n}- \sum_{n \leq x}d(n).
\end{eqnarray} 
Using the standard asymptotics for these finite sums, it yields
\begin{eqnarray}
\sum _{n\leq x} \left \{\frac{x}{n}\right \}&=&x \left(\log x + \gamma +O\left (\frac{1}{x} \right ) \right ) - \left(x\log x +(2 \gamma-1)x +O(x^{1/3}) \right ) \nonumber\\
&=&(1-\gamma)x+O\left (x^{1/3} \right )
\end{eqnarray} 
as claimed.
\end{proof}

The conjectured optimum error term is \(O\left(x^{1/4+\epsilon }\right)$, with $\epsilon >0\) any arbitrary small number. This coincides with the best possible error term required in the divisor problem and the circle problem, see \cite{CS07}. 
\\

\begin{lem} \label{lem15020}
	Let $x  \geq 1$ be a large number. Then
	\begin{equation}
	\sum _{n\leq x}\left \{\frac{x}{n} \right \}=(1-\gamma)x+\Omega_{\pm}\left(x^{1/4}\right).
	\end{equation}
	
\end{lem}

\begin{thm} \label{thm15005} Let $x \geq 1$ be a number, and let $\{x\}=x-[x]$ be the fractional part function. Then
	\begin{equation}
	\sum_{n \leq x} \frac{1}{n} \left \{\frac{x}{n} \right \}=a_0 \log x +a_1+O(x^{-1/2}) \nonumber,
	\end{equation}
	where $a_0=1-\gamma$ and $a_1$ are constants.
\end{thm}

\begin{proof} Let $F(t)=\sum_{n \leq t} \{t/n\}$. Using Theorem \ref{thm15000} and partial summation yield
	\begin{eqnarray}
	\sum_{n \leq x}\frac{1}{n} \left \{\frac{x}{n} \right \}
	&=&\int_1^x\frac{1}{t}\,dF(t) \nonumber \\
	&=& \frac{F(x)}{x}-F(1)+\int_1^x\frac{F(t)}{t^2}dt \nonumber \\
	&=& \frac{(1-\gamma)x+O(x^{1/2})}{x}-F(1)\\
	&&+\int_1^x \left( \frac{(1-\gamma)t+O(t^{1/2})}{t^2} \right ) dt \nonumber \\
	&=& a_0 \log x +a_1+O(x^{-1/2}) \nonumber,
	\end{eqnarray}
	where $a_0=1-\gamma$ and $a_1$ are constants.
\end{proof}

\subsection{Twisted Finite Fractional Sums}
The next Lemmas give estimates for the twisted summatory fractional Mobius function.
\\

\begin{lem} \label{lem15400}  Let \(x\geq 1\) be a large number, let \(\mu (n)\) be the Mobius function, and let $\{x\}=x-[x]$  be the fractional part function. Then,
\begin{equation}
\sum _{n \leq  x} \mu(n) \left \{\frac{x}{n} \right \}=-1+O\left(\dfrac{x}{\log ^2 x}\right).
\end{equation}
\end{lem}
\begin{proof} Let $F(x)=[ x ]$ be the largest integer function, and let $G(x)=1$ in Lemma \ref{lem40001}. Next, replace the integer part-fractional part identity:
\begin{eqnarray}
1&=&\sum_{n\leq x} \mu(n)\left [\frac{x}{n}\right] \nonumber\\
&=& \sum_{n\leq x} \mu(n)\left (\frac{x}{n}-\left \{ \frac{x}{n}\right \} \right )\nonumber\\
&=&x \sum_{n\leq x} \frac{\mu(n)}{n}- \sum_{n \leq x}\mu(n) \left \{ \frac{x}{n}\right \}  \\
&=& O\left (\dfrac{x}{\log^{2} x}\right )-\sum_{n\leq x} \mu(n) \left \{\frac{x}{n}\right \}\nonumber,
\end{eqnarray}
where the last line follows from Lemma \ref{lem21040}. Now, solve for the fractional Mobius sum.  \end{proof}

Almost the same calculatrion appears in \cite[p.\ 590]{LE01}.\\

\begin{lem} \label{lem15401}  Let \(x\geq 1\) be a large number, let \(\mu (n)\) be the Mobius function, and let \(\{
z \}\) be the fractional part function. Then,
\begin{equation}
\sum_{n\leq x} \frac{\mu(n)}{n}\left \{\frac{x}{n}\right\}=O(1) .
\end{equation}
\end{lem}

\begin{proof} Let \(V(x)=\sum _{n \leq  x} \mu (n)\{x/n\}=-1+O\left(x \log^{-2} x\right)\), see Lemma \ref{lem15400}. The integral representation yields
\begin{eqnarray}
\sum_{n\leq x} \frac{\mu(n)}{n}\left \{\frac{x}{n}\right\}
&=& \int_{1}^{x}\frac{1}{t} d V(t) \nonumber\\
&=& \frac{V(x)}{x}-V(1)+\int_{1}^{x}\frac{V(t)}{t^2} d t\\
&=&O(1)\nonumber.
\end{eqnarray}

Note that the integral
\begin{equation}                         
\int_{1}^{x}\frac{V(t)}{t^2}d t=\int_{1}^{x}\frac{-1+O\left(t\log ^{-2} t\right)}{t^2}d t=O\left (\frac{1}{\log x} \right).
\end{equation}
This verifies the claim.  
\end{proof} 

Similar calculations as in Lemmas \ref{lem15400} and \ref{lem15401} are given in \cite[p.\ 248]{MV07}.\\

\subsection{Quasibalanced And Balanced Average Orders}
The average order of the \textit{quasibalanced} fractional part function $\psi(x)=\{x\}-1/2$ appears in several problems,  see \cite{WA27}, \cite{SR85}, \cite{FS17}, et alii. But the applications for \textit{balanced} fractional part function $\psi_0(x)=\{x\}-(1-\gamma)$ are very rare.\\

\begin{lem} \label{lem15007} {\normalfont(Quasibalanced Average) } Let $x \geq 1$ be a number, and let $\{x\}=x-[x]$ be the fractional part function. Then
	\begin{equation}
	\sum_{n \leq x} \frac{\{x/n\}-1/2}{n} =(1/2-\gamma) \log x +c_0+O \left (\dfrac{1}{x^{1/2}} \right ),
	\end{equation}
	where $c_0$ a is constant.
\end{lem}

\begin{proof} Applications of Theorem \ref{thm15005} and the asymptotic formula for the harmonic finite sum yield
	\begin{eqnarray}
	\sum_{n \leq x} \frac{\{x/n\}-1/2}{n}
	&=&\sum_{n \leq x}\frac{1}{n}  \left \{\frac{x}{n} \right \}-\dfrac{1}{2}\sum_{n \leq x}\frac{1}{n}   \\
	&=& \left( (1-\gamma) \log x +c_1+O \left (\dfrac{1}{x^{1/2}} \right )\right )-\dfrac{1}{2}\left ( \log x + \gamma +O \left (\dfrac{1}{x} \right ) \right ) \nonumber \\
	&=& (1/2-\gamma) \log x +c_2+O \left (\dfrac{1}{x^{1/2}} \right ) \nonumber,
	\end{eqnarray}
	where $c_1$ and $c_2$ are constants.
\end{proof}

\begin{lem} \label{lem15009} {\normalfont(Balanced Average) } Let $x \geq 1$ be a number, and let $\{x\}=x-[x]$ be the fractional part function. Then
	\begin{equation}
	\sum_{n \leq x} \frac{\{x/n\}-(1-\gamma)}{n} =c+O \left (\dfrac{1}{x^{1/2}} \right ),
	\end{equation}
	where $c$ is a constant.
\end{lem}

\begin{proof} Applications of Theorem \ref{thm15005} and the asymptotic formula for the harmonic finite sum yield
	\begin{eqnarray}
	\sum_{n \leq x} \frac{\{x/n\}-(1-\gamma)}{n}
	&=&\sum_{n \leq x}\frac{1}{n}  \left \{\frac{x}{n} \right \}-(1-\gamma)\sum_{n \leq x}\frac{1}{n}  \nonumber \\
	&=& \left( (1-\gamma) \log x +c_1+O \left (\dfrac{1}{x^{1/2}} \right )\right )\\
	&&-\left ((1-\gamma) \log x + \gamma +O \left (\dfrac{1}{x} \right ) \right ) \nonumber \\
	&=& c_0+O \left (\dfrac{1}{x^{1/2}} \right ) \nonumber,
	\end{eqnarray}
	where $c_0$ is a constant.
\end{proof}

\subsection{Variance}
The average and the the quasibalanced fractional function $\psi(x)=\{x/n\}-1/2$ is slightly different from the uniform random variable on the interval $[-1/2,1/2]$. Howerver, the variance almost the same as a uniform random variable on the interval $[-1/2,1/2]$. \\

\begin{lem} \label{lem15011} {\normalfont(QuasiBalanced Variance) } Let $x \geq 1$ be a number, and let $\{x\}=x-[x]$ be the fractional part function. Then
	\begin{equation}
	\sum_{n \leq x} \left (\{\alpha n\}-1/2 \right )^2 =\frac{1}{12}x+O \left (x^{\varepsilon} \right ),
	\end{equation}
	where $\varepsilon>0$ is arbitrarily small constant depending on the irrational $\alpha >0$.
\end{lem}

\begin{proof} This requires the Fourier series of the Bernoulli polynomial $\psi(x)^2=(\{x\}-1/2)^2-1/12$. Summing the Fourier series over the range of integers yields
	\begin{eqnarray}
	\sum_{n \leq x} \left (\{\alpha n\}-1/2 \right )^2
	&=&\sum_{n \leq x} \left (\frac{1}{12} +\frac{1}{2 \pi^2}
    \sum_{m \geq 1}\frac{e^{i2 \pi \alpha mn }}{m^2} \right ) \nonumber \\
&=& \frac{1}{12}x +\frac{1}{2 \pi^2} \sum_{m \geq 1}\frac{1}{m^2} 	
    \sum_{n \leq x}e^{i2 \pi \alpha mn }\\
	&=&\frac{1}{12}x +  O \left ( x^{\varepsilon} \right ). \nonumber \\
	\end{eqnarray}
The error term $O(x^{\varepsilon})$ depends on the irrational $\alpha \in \R$; numbers with unbounded partial quaotients have the largest error term.
\end{proof}

This proof was known quite sometimes ago, see \cite{HL20}.

\subsection{Comparison Of Dirichlet And Walfisz Results} \label{ssect3}
The estimate of the quasibalanced fractional sum
\begin{equation}\label{key}
\sum_{n \leq x} \frac{\{x/n\}-1/2}{n} =(1/2-\gamma) \log x +c+O \left (\dfrac{1}{x^{1/2}} \right ),
\end{equation}
see Lemma \ref{lem15007} for a proof, is based on Dirichlet or delaValle Poussin result in Theorem \ref{thm15000}. This result seems to contradict a well known result described below. \\

\begin{thm} \label{thm15013} {\normalfont(Walfisz) } Let $x \geq 1$ be a number, and let $\{x\}=x-[x]$ be the fractional part function. Then
\begin{equation} \label{15600}
\sum_{n \leq x} \frac{\{x/n\}-1/2}{n} =O \left (\dfrac{\log x}{ \log \log x} \right ).
\end{equation}
\end{thm}

 The detailed and lenghty proof is provided in \cite[pp.\ 72--78]{WA27}. It is based on the Fourier series 
\begin{equation}
\psi(x)=-\frac{1}{\pi}\sum_{m\geq1} \frac{\sin(2 \pi mx)}{m}
=\left \{\begin{array}{ll}
\{x\}-1/2     &x \not \in \mathbb{Z},\\
0           &x \in \mathbb{Z},\\
\end{array}
\right.
\end{equation}
of the function $\psi(x)=\{x\}-1/2$, and an exponential sum estimate similar to
\begin{equation}
\sum_{n \leq x} e^{i \pi mx/n}=o(x).
\end{equation}

About fourty years later, the estimate in (\ref{15600}) was improved to 
\begin{equation}
\sum_{n \leq x} \frac{\{x/n\}-1/2}{n}=-\frac{1}{\pi} \sum_{n \leq x} \dfrac{1}{n}\sum_{m\geq1} \frac{\sin(2 \pi mx/n)}{m} =O \left (\log^{2/3}x \right ).
\end{equation}
This version and the associated results are often quoted in the literature, see \cite{SR82}, \cite{SS15}, \cite[p.\ 46]{TG15}, \cite[Section 2]{FS17}, et alii.

\section{Inversion Identities } \label{sec4}
\begin{lem} \label{lem40001} {\normalfont (Mobius summatory inversion)} Let \(F,G:\mathbb{N}\longrightarrow \mathbb{C}\) be complex-valued arithmetic
	functions. Then
\begin{equation}
F(x)=\sum_{n\leq x}G(x/n)\quad \text{and} \quad G(x)=\sum_{n\leq x} \mu(n)F(x/n)
\end{equation}
are a Mobius inversion pair.
\end{lem}

\begin{proof} Refer to \cite[p.\  237]{HW08}, \cite[p.\ 36]{MV07}, \cite[p.\ 25]{RH94}, \cite[p.\ 62]{SH83}, \cite[p.\ 35]{TG15}, et alii.      \end{proof}

\section{Results For The Ratio $\varphi(n)/n$} \label{sec5}
The corresponding normalized summatory totient function $\varphi(n)/n$ has the well known asymptotic formula $\sum _{n \leq  x} \varphi (n)/n = 6 \pi^{-2}x+O(\log
x)$, confer \cite[p.\ 36]{MV07}, and \cite[p.\ 229]{MR08}. Some earlier works on this problem appear in \cite{MH88}, \cite{PC30}, \cite{ES51}, \cite{KW10}, and similar references. An improved error term for the normalized summatory totient function is considered first.

\begin{thm}  \label{thm5.1} For large number \(x\geq 1,\) the average order for the normalized Euler totient function \(\varphi
	(n)/n\) has the asymptotic formula
\begin{equation}
\sum_{n\leq x}\frac{\varphi(n)}{n}=\frac{6}{\pi^2}x+O(1) .
\end{equation}
\end{thm}

\begin{proof} The analysis proceeds as usual, but improves on the last steps:
\begin{eqnarray}
\sum_{n\leq x} \frac{\varphi(n)}{n}
&=& \sum_{n\leq x}\sum_{d|n} \frac{\mu(n)}{d}\nonumber\\
&=&\sum_{d\leq x}\frac{\mu(n)}{d}\sum_{n\leq x,\,d|n} 1\\
&=& \sum_{d\leq x}\frac{\mu(n)}{d}\left [\frac{x}{d}\right]\nonumber, 
\end{eqnarray}
where \(\mu (n)\in \{ -1, 0, 1 \}\) is the Mobius function. Substituting the integer part/fractional part functions identity leads to
\begin{eqnarray}
\sum_{n\leq x} \frac{\varphi(n)}{n}
&=& \sum_{d\leq x} \frac{\mu(d)}{d} \left (\frac{x}{d}-\left \{\frac{x}{d}\right\} \right ) \\
&=& x\sum_{d\leq x} \frac{\mu(d)}{d^2}-\sum_{d\leq x} \frac{\mu(d)}{d} \left \{\frac{x}{d}\right\} \nonumber.
\end{eqnarray}
Now, using Lemmas \ref{lem21040} and \ref{lem15401} yields
\begin{eqnarray}
\sum _{n \leq  x} \frac{\varphi (n)}{n}&=&x\left(\frac{6}{\pi ^2}+O\left(\frac{1}{x \log ^2 x}\right)\right)+O(1) \nonumber\\
&=&\frac{6}{\pi ^2}x+O(1) .
\end{eqnarray}
This proves the claim.  \end{proof}

The standard proof for the average order \(\sum _{n \leq  x} \varphi (n)= 3\pi^{-2}x^2+O(x
\log  x)\) of \(\varphi (n)\) are due to Mertens, \cite{MF74}, \cite[p. 591]{LE01}. Currently, it is claimed that $\sum_{n \leq x}=3 \pi^{-2}x^2+O\left(x(\log
x)^{2/3}(\log \log x)^{4/3}\right)$, see \cite[p.\ 99?]{WA63}, \cite[p.\ 36]{MV07}, \cite[p.\ 46]{TG15}, and other authors.\\

\begin{proof} (Theorem \ref{thm1.1}) By Theorem \ref{thm5.1}, \(W(x)=\sum _{n \leq  x} \varphi (n)/n=6\pi ^{-2}x+O(1)\), and summation by part yields
\begin{eqnarray}
\sum_{n\leq x}\varphi(n)&=&\sum_{n\leq x} n\cdot\frac{\varphi(n)}{n}\nonumber\\
&=&\int_{1}^{x}t\: d W(t)\nonumber\\
&=& xW(x)+O(1)-\int_{1}^{x} W(t)d t\\
&=& x \left (\frac{6}{\pi^2}x+O(1) \right)-\int_{1}^{x} \left (\frac{6}{\pi^2}t+O(1) \right)d t \nonumber\\
&=&\frac{3}{\pi^2}x^{2}+O(x)\nonumber .
\end{eqnarray}
Quod erat demonstrandum.      
\end{proof}

\section{Results For The Ratio $n/\varphi(n)$} \label{sec6}
This section continues with the analysis of the error term of the average order for the reciprocal $1/\varphi(n)$ of the Euler totient function $\varphi(n)$. It proves  that the best error term is the same as that determined by Landau over a century ago, in \cite[p.\ 184]{LE00}. The simpler analysis for the ratio $n/\varphi(n)$ is considered first.

\begin{thm} \label{thm90001}
	Let $x \geq 1$ be a large number. Then, the average order of the ratio $n/\varphi(n)$ has the asymptotic formula	
	\begin{equation}	
	\sum_{n \leq x} \frac{n}{\varphi(n)} =a_0x+O(\log x),
	\end{equation}
	where $a_0=\zeta(2)\zeta(3)/\zeta(6)$ is a constant.
\end{thm}

\begin{proof} The result is derived using the identity $\sum_{d \mid n}\mu^2(d)/\varphi(d) $. Substituting this formula, and reversing the order of summation yield
	
	\begin{eqnarray}	
	\sum_{n \leq x} \frac{n}{\varphi(n)}  &=&  \sum_{n \leq x} \sum_{d \mid n}\frac{\mu^2(d)}{\varphi(d)}  \nonumber \\
	&=& \sum_{d \leq x}\frac{\mu^2(d)}{\varphi(d)} \sum_{n \leq x, \;d \,\mid\, n}1\\
	&=& \sum_{d \leq x}\frac{\mu^2(d)}{\varphi(d)} \left ( \frac{x}{d}-\left \{ \frac{x}{d}\right \} \right ) \nonumber \\
	&=& x\sum_{d \leq x}\frac{\mu^2(d)}{d\varphi(d)} - \sum_{d \leq x}\frac{\mu^2(d)}{\varphi(d)}  \left \{ \frac{x}{d} \right \}  \nonumber .	\end{eqnarray}	
	The first finite sum 
	\begin{eqnarray}	
	x\sum_{n \leq x}\frac{\mu^2(n)}{n\varphi(n)}&=& x\left ( \sum_{n\geq 1}\frac{\mu^2(n)}{n\varphi(n)}-\sum_{n > x}\frac{\mu^2(n)}{n\varphi(n)}\right )  \nonumber \\ &=&c_0x+O\left ( \frac{1}{x} \right ).	
	\end{eqnarray}
	The constant $a_0>0$ is expressable in terms of zeta functions as	
	\begin{equation}	
	\sum_{n \geq 1}\frac{\mu^2(n)}{n\varphi(n)}  = \prod_{p \geq 2} \left ( 1+\frac{1}{p(p-1)} \right )  =\frac{\zeta(2)\zeta(3)}{\zeta(6)}.	
	\end{equation}	
	The second finite sum
	\begin{equation} \label{90011}	
	\sum_{n \leq x}\frac{\mu^2(n)}{\varphi(n)}  \left \{ \frac{x}{n} \right \}  \gg \log x	
	\end{equation}	
	is alway positive and exhibits no cancellations.
	
\end{proof}

The form of the error term in (\ref{90011}) concretely proves that it cannot be improved, see also Subsection \ref{ssect3} for a related discussion. \\

\begin{thm} \label{thm90007}  {\normalfont (Same as Theorem 1.3)}  Let $x \geq 1$ be a large number. Then, the average order of the ratio $1/\varphi(n)$ is as follows.
	
	\begin{equation} \label{90013}	
	\sum_{n \leq x}  \frac{1}{\varphi(n)} =c_0+c_1\log x +O\left (\frac{\log x}{x} \right ),
	\end{equation}	
	where
	$c_0$ and $c_1$ are constants.	
\end{thm}

\begin{proof}   The result is derived Theorem \ref{thm90001} by partial summation. More precisely, let $R(t)=\sum_{n \leq x} n/\varphi(n)=a_0x+O(\log x)$. Then
	\begin{eqnarray}	
	\sum_{n \leq x} \frac{1}{\varphi(n)}  &=&\sum_{n \leq x}\frac{1}{n} \frac{n}{\varphi(n)}  \nonumber \\
	&=&  \int_1^x \frac{1}{t} d R(t)  \\
	&=& \frac{R(t)}{t} \, \Biggr |_1^x+ \int_1^x \frac{R(t)}{t^{2}} d t ,\nonumber \\
	&=& \frac{a_0x+O\left (\log x \right )}{x} +a_1 + \int_1^x \frac{a_0t+O\left (\log t \right )}{t^{2}} d t , \nonumber \\
	&=& c_0+c_1\log x +O\left (\frac{\log x}{x} \right ) \nonumber,
	\end{eqnarray}	
	where
	$a_0=\zeta(2)\zeta(3)/\zeta(6),a_1=-R(1),c_0,$ and $c_1 $ are constants.

\end{proof}

The work in \cite{SR82} is devoted to improving the error term from $O((\log x)/x)$ to $O((\log x)^{2/3}/x)$. This analysis was based on the estimate
\begin{equation}
\sum_{n \leq x} \frac{\{x/n\}-1/2}{n}=O \left ((\log x)^{2/3} \right ),
\end{equation}
Refer to subsection \ref{ssect3} for a discussion on this estimate. However, by Theorem \ref{thm90007}, this seems to be impossible since the error term in (\ref{90013}) satisfies $ \gg (\log x)/x$, confer (\ref{90011}).

\section{Abridged History Of The Error Term} \label{sec7}
Recall that by definition, the error term is given by 
\begin{equation}
E(x)=\sum_{n \leq x} \varphi(n)-\frac{3}{\pi^{2}}x^2
\end{equation}. 

The earliest estimate of the error term $E(x)=O(x^{1+\varepsilon})$ was computed by Dirichlet, followed by Mertens as $E(x)=O(x\log x)$, see \cite{MF74}, and later $E(x)=O(x(\log x)^{2/3}(\log \log x)^{4/3})$ was computed by Walfisz, see \cite[p.\ 99]{WA63}. The assertions that $E(x)\ne o(x \log \log \log x)$, and $E(x)= \Omega(x \log \log \log \log x)$, appear in \cite{PC30}, and \cite{ES51} respectively. Moreover, there is a conjecture that $R(x)=O(x \log \log x)$, and the omega estimate
\begin{equation}
\sum_{n\leq x}\varphi(n)-\frac{3}{\pi^2}x^2=\Omega_{\pm}(x \sqrt{\log \log x})
\end{equation} 
confer \cite{MH88}. The analysis given in those papers are for the error terms over a short interval, namely, $R(x)=\sum_{x\leq n \leq x+y} \varphi(n)$, where $y=O(\log \log x)$, not for $E(x)=\sum_{n \leq x} \varphi(n)-3\pi^{-1}x^2$.\\ 

Futhermore, are a few other recent result such as the smoothed omega estimate
\begin{equation}
\sum_{n\leq x}\varphi(n)\log \left(\frac{x}{n} \right )-\frac{3}{\pi^2}x^2=\Omega_{\pm}(x^{1/2} \log \log \log x)
\end{equation}
was proved in \cite{KW10}. 

\newpage\section{Problems}
1. Find the exact value of the main term of the finite sum
\begin{equation} \label{92011}	
	\sum_{n \leq x}\frac{\mu^2(n)}{\varphi(n)}  \left \{ \frac{x}{n} \right \} 	
	\end{equation}	

2. Let $\psi(x)=\{x\}-1/2$. Comupte the Fourier series
	\begin{equation}
	\psi(x)=\frac{-1}{ \pi} \sum_{n \geq 1}\frac{\sin( 2\pi x n)}{n}.
	\end{equation} 

3. Let $\psi(x)=\{x\}-1/2$. Comupte the Fourier series
	\begin{equation}
	\psi(x)^2=\frac{1}{12} +\frac{1}{2 \pi^2} \sum_{n \geq 1}\frac{e^{i 2 \pi x n}}{n^2}.
	\end{equation}

4. Let $\alpha \in \R- \Z$ be irrational, with unbounded partial qutients, estimate the error term $\sum_{n \leq x} \left (\{\alpha n\}-1/2 \right )^2 -x/12 =O \left (x^{\varepsilon} \right )$ in Lemma \ref{lem15011}.

\end{document}